\documentclass[reqno,11pt,a4paper]{amsart}
\usepackage[margin=2.5cm]{geometry}
\usepackage{pgfplots}

\usepackage{hyperref}
\usepackage{pdflscape}
\usepackage{longtable}
\usepackage{booktabs}
\usepackage{colortbl}
\usepackage{arydshln}
\usepackage{verbatim}
\usepackage{enumitem}
\usepackage{multirow}
\usepackage{makecell}
\usepackage{pdfsync}

\usepackage{amsthm}
\usepackage{hyperref}
\usepackage{enumerate}
\usepackage{url}
\usepackage{pdflscape}
\usepackage{longtable}
\usepackage{booktabs}
\usepackage{colortbl}
\newcommand{\evnrow}{\rowcolor[gray]{0.95}}
\newcommand{\oddrow}{}

\usepackage[margin=2.5cm]{geometry}

\usepackage{pdflscape}
\usepackage{longtable}
\usepackage{booktabs}
\usepackage{colortbl}
\usepackage{arydshln}
\usepackage{verbatim}
\usepackage{enumitem}
\usepackage{mathdots}

\usepackage[english]{babel}
\usepackage{lipsum}
\usepackage{amsthm}
\usepackage{todonotes}
\usepackage[vcentermath]{youngtab}

\theoremstyle{plain}
\newtheorem{lema}{Lemma}[section]
\newtheorem{prop}[lema]{Proposition}
\newtheorem{thrm}[lema]{Theorem}

\theoremstyle{remark}

\newtheorem{rmk}[lema]{Remark}

\theoremstyle{definition}
\newtheorem{dfn}[lema]{Definition}

\def\ZZ{{\mathbb Z}}

\def\CC{{\mathbb C}}

\def\AA{{\mathbb A}}

\def\PP{{\mathbb P}}
\def\wgr{\mathrm{wGr(2,5)}}
\def\Gr{\mathrm{Gr}}
\def\wp2{\mathrm{w}(\PP^2\times \PP^2)}

\def\wG{\mathrm{w}\mathcal{G}}
\def\wP{\mathrm{w}\mathcal{P}}

\def\w{{\mathrm{w}}}

 \newcommand{\into}{\hookrightarrow}

 \newcommand{\Oh}{\mathcal O}

 \newcommand{\Si}{\Sigma}

\newcommand{\PxP}{\PP^2 \times \PP^2}
\usepackage{amstext}
\usepackage{amssymb}

\begin{document}

%--------------------------------------------------------------------------------------------

\author[M.I.~Qureshi]{Muhammad Imran Qureshi}
\address{Deptartment of Mathematics \& Statistics, King Fahd University of Petroleum and Minerals, Saudi Arabia}
\email{imran.qureshi@kfupm.edu.sa}
%-------------------------------------------------------------------------------
%--------------------------------------------------------------------------------------------
\title{Smooth Fano $4$-folds in Gorenstein formats}
%--------------------------------------------------------------------------------------------
\begin{abstract}
We construct some new deformation families of four dimensional Fano manifolds of index $1$ in some known classes of Gorenstein formats.  These families have explicit descriptions in terms of equations, defining  their image under the anti-canonical embedding in some weighted projective space. The constructed families  have relatively smaller anti-canonical degree than most other known families of smooth Fano 4-folds.

\end{abstract}
%--------------------------------------------------------------------------------------------
\maketitle
%\tableofcontents
%--------------------------------------------------------------------------------------------

\section{Introduction}

A  projective  algebraic variety \(X\) with an ample anti-canonical divisor class \(-K_X\) is called a Fano variety. The Fano index \(i\) is the largest integer such that \(-K_{X}=iD\), for some ample divisor \(D\) on \(X\). Fano varieties are one of the central topics of research in algebraic geometry in general and in  classification problems  in particular. It is well known that there are only finitely many deformation families of smooth Fano varieties in each dimension \cite{KMM}. In dimension  less than or equal to three the classification has been completed \cite{Isk1,Isk2,Isk-Pro,MM,MM2}. There are 1, 10 and 105 deformation families of smooth Fano varieties in dimension 1, 2 and 3 respectively.  

In dimension greater than or equal to 4 the full classification is still to be completed. The  complete classification of smooth Fano 4-folds of  index  greater than or equal to 2 is known and there are $35$ deformation families of such Fano 4-folds \cite{Fuj1,Fuj2,Fuj3,Fuj4,Isk1,Isk-Pro,Kob1,Wilson,Wis}, listed in \cite{Q-periods}.   The index 1 case is still not complete, although there are  a number of partial classification results. The toric Fano 4-folds have been classified by Batyrev \cite{Bat-toric-4} and found 123 deformation families. 
One of the larger set of examples was constructed by     Coates, Kasprzyk and Prince: in \cite{CKP}  they constructed $527$ new deformation families of Fano 4-folds as complete intersections in toric varieties and one more in \cite{CKP-LI} by using the Laurent inversion. Another collection of    $141$ deformation families have been given  by Coates, Kasprzyk and Kalashnikov  \cite{Klash} as quiver flag zero loci in quiver flag varieties. The  smooth Fano 4-folds of Picard rank 2 with hypersurface Cox ring have been classified by Hausen, Laface and Mauz  \cite{HLM}: finding     $  17$ new deformation families of smooth Fano 4-folds of index 1. Some of the earlier index 1 examples have been constructed by K\"uchle \cite{Kuchle-95,Kuchle-97} as  sections of homogeneous vector bundles over Grassmannians and complete intersections in weighted projective spaces.

In this article, we aim to contribute to the classification of smooth Fano 4-folds of index 1, which can be thought of as an  extension of  K\"ucle's lists as  some  of his examples also appear in our list. We construct some new deformation families of smooth Fano 4-folds of index 1 as weighted complete intersections of some known classes of Gorenstein formats (Definition \ref{dfn:GF}). As a starting point,  we use a computer algebra system   to search for candidate families of smooth Fano 4-folds by using the algorithmic approach  \cite{QJSC,BKZ}. Then we prove the existence of these 4-folds by analysing the explicit equations of these candidate varieties.

 %--------------------------------------------------------------------------------------------
\subsection*{Summary of results}
\label{Smooth-Fano4}
In total, we obtained   25 candidate families of smooth Fano 4-folds  by using the computer search routine \cite{QJSC, BKZ}. Four as hypersurfaces,  8 as codimension 2 complete intersections, one each complete intersection in codimension 3 and 4, four in Gr(2,5) format, four in $\PxP$ format and three in codimension 4 \(\Gr(2,5) \cap \mathcal H\) format.  Among these one candidate codimension 2 complete intersection and one  candidate examples in  \(\Gr(2,5) \cap \mathcal H\) format failed to be smooth. All the hypersurfaces and complete intersection examples already appeared in  \cite{Kuchle-97}, so new deformation families  appear as non-complete intersection Fano 4-folds. 
\begin{thrm} \label{main}There exist at least ten  families of smooth Fano 4-folds of index 1 such that their images under the anti-canonical  embedded in a weighted projective space  can be described    as non-complete intersection varieties, given in Table \ref{Fano4-Formats}. In four of these cases they can be  be described in \(\PP^7(w_i)\) by using \(\Gr(2,5)\) format, in two cases by  \(\Gr(2,5)\cap \mathcal H   \) format in \(\PP^8(w_i)\) and in four cases by using \(\PxP\) format in \(\PP^8(w_i)\). The families 1--6 have Picard rank 1 and 7--10 have Picard rank 2. 
\end{thrm}
\begin{longtable}{>{\hspace{0.5em}}llccccr<{\hspace{0.5em}}}

\caption{ Smooth Fano 4-folds in Gorenstein formats} \label{Fano4-Formats}\\
\toprule
\multicolumn{1}{c}{No}&\multicolumn{1}{c}{Format}&\multicolumn{1}{c}{ $(-K_X)^4$}&\multicolumn{1}{c}{$h^0(-K_X)$}&\multicolumn{1}{c}{Eq Degs,  $\w\PP$}& \multicolumn{1}{l}{Weight Matrix}\\
\cmidrule(lr){1-1}\cmidrule(lr){2-2}\cmidrule(lr){3-4}\cmidrule(lr){5-6}
\endfirsthead
\multicolumn{7}{l}{\vspace{-0.25em}\scriptsize\emph{\tablename\ \thetable{} continued from previous page}}\\
\midrule
\endhead
\multicolumn{7}{r}{\scriptsize\emph{Continued on next page}}\\
\endfoot
\bottomrule
\endlastfoot

\evnrow $ 1$&$\Gr(2,5) $& 13&8&$ \begin{array}{@{}l@{}}\oddrow X_{2,3^4}\\\quad\subset\PP(1^8)\end{array}$&$\footnotesize\begin{matrix} 1&1&1&2\\ &1&1&2\\ &&1&2\\&&&2 \end{matrix}$\\
\oddrow $ 2$& &10&7& $ \begin{array}{@{}l@{}} X_{3^4,4}\\\quad\subset\PP(1^7,2)\end{array}$&$\footnotesize\begin{matrix} 1&1&1&1\\ &2&2&2\\ &&2&2\\&&&2 \end{matrix}$\\
\evnrow $ 3$& &7&$6  $& $ \begin{array}{@{}l@{}}\oddrow X_{2,3^4}\\\quad\subset\PP(1^6,2^2)\end{array}$&$\footnotesize\begin{matrix} 1&1&2&2\\ &1&2&2\\ &&2&2\\&&&3 \end{matrix}$\\
\oddrow $ 4$& &5&$ 5 $&$ \begin{array}{@{}l@{}} X_{4^5}\\\quad\subset\PP(1^5,2^3)\end{array}$&$\footnotesize\begin{matrix} 2&2&2&2\\ &2&2&2\\ &&2&2\\&&&2 \end{matrix}$\\
\evnrow $ 5$&$\Gr(2,5)\cap \mathcal H$ &15&$ 9  $&$ \begin{array}{@{}l@{}}\evnrow X_{2^5,(3)}\\\quad\subset\PP(1^9)\end{array}$&$\footnotesize\begin{matrix} 1&1&1&1\\ &1&1&1\\ &&1&1\\&&&1 \end{matrix}$\\
\oddrow $ 6$& &10&$8   $&$ \begin{array}{@{}l@{}}\oddrow X_{2^5,(4)}\\\quad\subset\PP(1^8,2)\end{array}$&$\footnotesize\begin{matrix} 1&1&1&1\\ &1&1&1\\ &&1&1\\&&&1 \end{matrix}$\\
\evnrow  $ 7$&$ \PxP $&$17$ &9&$ \begin{array}{@{}l@{}} X_{2^3,3^6}\\\quad\subset\PP(1^9) \end{array}$&$\footnotesize\begin{matrix} 1&1&1\\ 1&1&1\\ 2&2&2 \end{matrix}$\\
\oddrow $ 8$& &11&7&$ \begin{array}{@{}l@{}}\oddrow X_{3^6,4^3}\\\quad\subset\PP(1^7,2^2)\end{array}$&$\footnotesize\begin{matrix} 1&1&1\\ 2&2&2\\ 2&2&2 \end{matrix}$\\
\evnrow $ 9$&& $10$&7&$ \begin{array}{@{}l@{}} X_{2,3^4,4^4}\\\quad\subset\PP(1^7,2^2)\end{array}$&$\footnotesize\begin{matrix} 1&1&2\\ 1&1&2\\ 2&2&3 \end{matrix}$\\
\oddrow $ 10$& & $5$&5&$ \begin{array}{@{}l@{}}\oddrow X_{4^9}\\\quad\subset\PP(1^5,2^4)\end{array}$&$\footnotesize\begin{matrix} 2&2&2\\ 2&2&2\\ 2&2&2 \end{matrix}$\\

 \end{longtable}
 
 This list of examples in not a formal complete classification of Fano 4-folds in these Gorenstein formats but it is very unlikely that we  get any further such examples.  In fact, we searched for about half a million candidate embeddings of Fano 4-folds by using computer algebra, combining  all cases. The examples  1 and 5 also appeared  in \cite{Kuchle-97} and the rest of the examples are new deformation families of smooth Fano 4-folds.

 It is well known that the plurigenera of  the Hilbert series $P_X(t)=\sum h^0(-nK_X)t^n  $ are deformation  invariant  \cite{Inv:pluri}. In particular, by using Reimann-Roch formula for smooth  4-folds  one can establish that   knowing\ the first plurigenus \(h^0(-K)\) and  the anti-canonical degree \((-K_X)^4\) is sufficient to distinguish between non deformation equivalent families, see Proposition \ref{distinguish}. For the known families of smooth Fano 4-folds the list of these invariants can be found in  \cite{Kuchle-95,Kuchle-97,Q-periods, CKP,CKP-LI,Klash,HLM}. We obtain at least 6 new deformation families of smooth Fano 4-folds of index 1.   The lists we obtain are conjecturally complete classifications of such varieties in these formats. 

 A linear section \(Y\) of each of the  Fano 4-folds in Table \ref{Fano4-Formats} is  a smooth Calabi--Yau 3-fold. Therefor by Lefschetz's hyperplane theorem, the  Picard rank of a Fano 4-fold \(X\) will be equal to that of  its Calabi--Yau 3-fold section \(Y\). The Picard rank for the corresponding \(Y \) have been calculated by using \cite[Theorem 2.5]{DNFF} and  computer algebra package \cite{Ilten} in Macaulay2.

% \footnote{Can this be proved that these families are indeed different. In Gr(2,5) case we know \(h^{1,1}=1\) by using \cite{BF} and can compute \(h^{1,3}\) using computer algebra in both \(\PxP\) and Gr(2,5) formats. We probably need a formal proof that $h^{1,1}=2 $ in \(\PxP\) examples but may follow from \cite{BF} or by using computer algebra on its linear section (Calabi--Yau 3-fold). So in principal we are left to calculate either \(e(X)\) or \(h^{2,2}(X)\) to prove this result. Their Calabi--Yau 3-fold linear section have different Hodge numbers so it is very likely they are not in the same deformation family.  } .  
     
     \begin{rmk}
     A pair of examples  \#2 and  \#9 and another pair of examples \#4 and \#10 have same numerical invariants \(h^0(-K_X)\) and \((-K_X)^4\) but they lie in different codimension. We  expect them to lie in  different deformation families;  a  similar phenomenon was observed for  some terminal Fano 3-folds   in \cite{BKQ}.   The corresponding Calabi--Yau 3-folds sections of  pairs \#2, \#9 and \#4, \#10 have distinct Hodge numbers \cite{MNQ}  and lie in different deformation families,   which provide an evidence that the  corresponding Fano 4-folds also belong to distinct deformation families. 
     \end{rmk}
\begin{rmk} Our computer search routine also recovered  the 13 examples of  smooth Fano 4-folds of index 1 which are  hypersurfaces or  complete intersections in weighted projective spaces.  However, we do not list them as a part of a theorem since they are very well known examples that appeared in \cite{Kuchle-97} and further studied in \cite{PS-bounds}. We list them for the reader in  Table \ref{Fano4-CI}.    
\begin{table} 
\caption{Smooth Fano 4-fold hypersurfaces and complete intersections list of  \cite{Kuchle-95}}
\vspace{-2mm}
\label{Fano4-CI} \[
\renewcommand{\arraystretch}{1.2}\begin{array}{|c|c|c|c|l|}
\hline
\text{Format} & \text{codimension} & h^0(-K_X)  & (-K_X)^4 & \text{Fano 4-fold}  \\\hline  \multirow{4}{*}{Hyp. Surf.} & \multirow{4}{*}{1} &{6} &5 & X_{5}\subset \PP(1^6)  \\\cline{3-5}
    & & \multirow{2}{*}{5} &3 & X_{6}\subset \PP(1^5,2) \\\cline{4-5}
    & &  & 2& X_{8}\subset \PP(1^5,4) \\\cline{3-5}
       & & 4 & 1& X_{10}\subset \PP(1^4,2,5) \\\cline{1-5}

 \multirow{11}{*}{Comp. Int.} & \multirow{7}{*}{2} &\multirow{2}{*}{7} &9 & X_{3^2}\subset \PP(1^7)  \\\cline{4-5}
    & &  &8 & X_{4,2}\subset \PP(1^7) \\\cline{3-5}
    & & 6 & 6& X_{4,3}\subset \PP(1^6,2) \\\cline{3-5}
& &5 &4 & X_{4^2}\subset \PP(1^5,2^2)\\\cline{3-5}
& &6  & 4&X_{6,2}\subset \PP(1^6,3) \\\cline{3-5}
& & 4 &2& X_{6,4}\subset \PP(1^4,2^2,3) \\\cline{3-5}
& & 3 & 1&X_{6^2}\subset \PP(1^3,2^2,3^2)\\
\cline{2-5}
& 3& 8 & 12&X_{2^2,3}\subset \PP(1^8)\\\cline{2-5}
& 4& 9 & 16&X_{2^4}\subset \PP(1^9)\\\cline{1-5}
 \hline

\end{array}
\]
\end{table}
\end{rmk}
\begin{rmk} We did search for examples of smooth Fano 4-folds in some other well known  classes of Gorenstein formats, namely in  Grassmannians Gr(2,6) \cite{QS}, Lagrangian Grassmannians LGr(3,6)\cite{QS2}, a two step flag variety in \( \CC^{4}\) \cite{QS2}, Othogonal Grassmannians OGr(5,10)\cite{wg} and weighted homogeneous \(F_4\) variety   \cite{QJGP} but no new candidate examples  of smooth Fano 4-folds were found.  
\end{rmk}
%-------------------------------------------------%%
\section{Definitions and notations}
%-------------------------------------------------%%
A weighted projective variety  \(X \into \PP^N(w_i)\) of codimension \(c\) is called \textit{wellformed} if it does not contain a singular stratum of   codimension \(c+1\). 
 A format is roughly a way of representing the equations of varieties. For example, the Segre embedding of \(\PxP\)  can be described as \(2 \times 2 \) minors of the size 3 matrix. A more formal definition of Gorenstein format  is given below. 
\begin{dfn}\cite{BKZ}\label{dfn:GF} A codimension \(c\) \emph{Gorenstein format}  \(\mathcal F\) is a triple \(\left(\widetilde{V } 
 , \mathcal{R}, \mu \right)\) which consists of a codimension \(c\) affine Gorenstein variety \(\widetilde V\subset \AA^n\), a minimal graded free resolution \(\mathcal R\) of \(\Oh_{\widetilde V}\) as a graded \(\Oh_{\AA^n}\) module, and a \(\CC^*\)-action \(\mu\) of strictly positive weights on \(\widetilde V \).  
 \end{dfn}
 We only consider those Gorenstein formats where the action \(\mu\) leaves the variety \(\widetilde V\) invariant  and the  free resolution \(\mathcal R\) is equivariant for the action. The varieties defined below as    \ref{dfn:WG} and \ref{dfn:P2xP2} are examples of such Gorenstein formats.
 
\begin{dfn}\cite{wg}\label{dfn:WG} Consider the  Pl\"ucker embedding 
  \(\Gr(2,5) \into\PP^9\left(\bigwedge^2 \CC^5\right) \) of Grassmannians of 2-planes in \(\CC^5\).  For a choice of  vector of  $w:=(a_{1},\cdots,a_5)$  where all \(a_i\in \frac12\ZZ\) satisfying   $$a_i+a_j>0,\; 1\le i<j\le 5,$$ one can define the weighted Grassmannian \(\w\Gr(2,5)\)  as quotient of affine cone minus the vertex     $\widetilde{\Gr(2,5)} \backslash\{\underline 0\} $ by $\CC^\times$ action given by: 
   $$\mu:x_{ij}\mapsto \mu^{a_i +a_j }x_{ij}.$$ \end{dfn}
\noindent Therefore we get the embedding
\begin{equation}\label{embed}\wgr\into \PP(\left\{w_{ij}: 1\le i <j\le 5, w_{ij}=a_i+a_j\right\}).  \end{equation} We will use \(\wG\) to denote the \(\w\Gr(2,5)\). The image of \(\wG\) under the embedding \eqref{embed} can be defined by  five maximal Pfaffians of $5\times 5$ skew symmetric matrix  
   \begin{equation}\label{eq:pf_mat}\left(\begin{matrix} 
x_{12}&x_{13}&x_{14}&x_{15}\\ 
&x_{23}&x_{24}&x_{25}\\ 
&&x_{34}&x_{35}\\ 
&&& x_{45} \end{matrix}\right)
,\end{equation}
where we omit writing down the diagonal of zeros and lower triangular part of the skew symmetric matrix.  
If   $\wG$ does not contain a 5 dimensional singular locus of \(\PP(w_i)\)  then the  canonical divisor  class is given by  \begin{equation}\label{eq:K_G}K_{\wG}= \left(-\frac12\sum_{1\le i<j\le 5} w_{ij} \right)H,\end{equation}
for an ample divisor \(H\). 
Another variant of  this format is the  \(\Gr(2,5) \cap \mathcal H\) format which is  a non-quasilinear hypersurface section of   $\Gr(2,5)$ format.   
\begin{dfn}\cite{BKQ,Sz05} \label{dfn:P2xP2} Let \(\Si\) be the Segre embedding of  \[\PxP\into \PP^8(x_{ij}),\;\; 1\le i,j\le 3\]  and  
 consider a pair of  half integer vectors   $b=(b_1,b_2,b_3)$ and $c=(c_1,c_2,c_3)$ 
$$b_i+c_j > 0,b_i\le b_j \text{ and } c_i\le c_j \text{ for } 1\le i\le j \le 3. $$  
 Then   quotient  of the punctured affine cone $\widetilde{\Si}\backslash\{\underline 0\}$  by  $\CC^\times$:
$$\mu:x_{ij}\mapsto\mu^{b_i+c_j}x_{ij}, \; 1\le i,j\le 3,$$ is called a  weighted \(\PxP\) variety, which we will denote by \(\wP\).  Thus for   a  parameter $p=(b,c)=(b_1,b_2,b_3;c_1,c_2,c_3), $ we get the embedding $$\wP\into \PP^8\left(w_{ij}: w_{ij}=b_i+c_j;1\le i,j\le 3\right). $$ 
 \end{dfn}
 The equations of image are well known  $2\times 2$ minors of the $3\times 3$ matrix, which we usually refer to as  weight matrix and write it as
\begin{equation}\label{eq:wtmx_P2}
\begin{pmatrix} w_{11} & w_{12} & w_{13} \\ w_{21} & w_{22} & w_{23} \\ w_{31} & w_{32} & w_{33} \end{pmatrix} \text{ where }w_{ij}=b_i+c_j; 1\le i,j\le 3.
\end{equation}
 If $\wP$ is wellformed then the  canonical divisor class is given by \begin{equation}\label{eq:K_P2}K_{\wP}=\left(-\sum_{ i=1}^3 w_{ii} \right)H,\end{equation}
for an ample divisor \(H\).

% \subsection*{Notation and Conventions}
% \begin{itemize}
% \item We work over the field of complex numbers \(\CC\). 
% \item We use the same notation for canonical divisor class \(K_X\) and canonical sheaf \(\om_X\), if no confusion can arise. We usually write \(K_{X}=\Oh(n)\) to represents \(K_X = nD\).   
% \end{itemize}
%%%%%%%%%%%%%%%%%%%%%%%%%%%%%%%%%%%%%%%%%%%%%%%%%%%%%%%%%%%%%%%%
\section{families of Fano smooth 4-folds}
In this section we provide a proof of the  theorem \ref{main} by providing details of the calculations in three cases. The rest of the cases can be computed by performing similar calculations. \\
\subsection{Example \#3}
 For \(w=\frac12(1,1,1,3,3)\) we get the embedding  \(\wG\into \PP\left(1^3,2^6,3\right)\). Let \(x_{1},x_2,x_3\) be weight one variables, \(y_1,\cdots,y_6\) be weight 2 variables and \(z\) be the  weight 3 variable.  Now it is evident  that \(\wG\) does not contain any 5-dimensional singular locus of the ambient \(\PP^9(w_i)\):  the weight 3 locus is just an orbifold point and weight \(2\) locus  describes a cubic 3-fold in \(\PP^5\), so it is wellformed.  Thus by using formula \eqref{eq:K_G} we have \(K_{\wG}=\Oh(-9)\).  
 Let   \(Y_1\subset \PP\left(1^6,2^6,3\right)\) be a projective cone over \(\wG\) with vertex \(\PP^2\), i.e. we add three variables of weight 1 to the ambient \(\PP^9(w_i)\) which are not involved in any defining equations of \(\wG\).  \(Y_1\) is a 9-dimensional variety with \(K_{Y_1}=\Oh\left(-12\right)\). 

\noindent We take a complete intersection of \(Y_1\) with \(4\) general quadrics to get a Fano 5-fold \[Y_2\subset \PP\left(1^6,2^2,3\right) \text{with } K_{Y_2}=\Oh(-12+2\times 4)=\Oh(-4) ,\] by using adjunction formula.  Now the base locus of the linear system of quadrics \(|\Oh(2)|\) contain a single point which is a coordinate point of the variable \(z\). Moreover \(Y_2\cap \PP(2,2)\) is empty locus and only  singular point  on \(Y_2\) is  \(\frac{1}{3}(1,1,1,2,2)\) quotient singularity. As a last step we take an  intersection of \(Y_{2}\) with a general  cubic to get a Fano 4-fold  
\[X\subset \PP(1^6,2^2) \textrm{ with } K_X=\Oh(-1).\]
As \(X\) does not contain any singular point of the \(\PP(1^6,2^2)\) so \(X\) is a smooth Fano 4-fold of index 1. By using the Hilbert series one can easily show that  \((-K_X)^4=7,\) and \(h^0(-K_X)=6 \) is evident from the embedding. 
 
\subsection{Example \#6} The Grassmannian \(\Gr(2,5)\) has the embedding in \(\PP^9(x_1,\cdots,x_{10})\). It is a 6-fold with \(K_X=\Oh(-5)\). Let \(Y_1\) be variety obtained by taking a
  cone of  weight 2 over it, i.e. we have the embedding 
  \[Y_1=\mathcal C^2\Gr(2,5)\into \PP^{10}(1^9,2)\] where the new variable \(y\) of weight 2 is not involved in any defining equations. \(Y_1\) is a singular 7-fold with \(K_{Y_1}=\Oh(-7)\). Now we take a general quartic section \[Q_4=y^2+f_4(x_i,y), 1\le i\le 10\] of \(Y_1\)  to get a 6-fold \[Y_2\subset \PP(1^{10},2) \textrm{ with } K_{Y_2}=\Oh(-3).\]
The linear system \(|\Oh(4)|\) has empty base locus. The 6-fold \(Y_2\) is a codimension 4 smooth variety; since  the weight 2 points have been removed by the quartic section. Now we take two hyperplane sections of \(Y_2\) to get a smooth Fano 4-fold \(X\subset \PP(1^8,2)\) of index \(1\)  with \(h^0(-K_X)=8\) and \(-K_X^4=10\). 

\subsection{Example\#10} For a choice of parameter \(w=(1,1,1;1,1,1) \) we have \(\wP \into \PP\left(2^9\right) \) which is a priory a non well-formed 4-fold. Consider a projective cone \(Y_{1}\) over \(\wP\) with vertex \(\PP^4\) then we have a 9-fold 
\[Y_1\subset \PP(1^5,2^9) \text{ with  } K_{Y_1}=\Oh(-11).\]
The 9-fold $Y_1$ is wellformed, though it contains the orbifold locus of dimension 4 defined by weight 2 variables and further  \(\PP^4\) given by cone variables. Then the  quasilinear section of  \(Y_1\) with 5 general quadrics is a Fano 4-fold 
\[X=Y_1\cap \{ \displaystyle\cap_{i=1}^5 Q_i\}\subset \PP(1^5,2^4) \text{ with } K_{X}=\Oh(-11+10)=\Oh_X(-1).\]
i.e.  \(X\) is a  Fano 4-fold of index 1. The intersection  \(X\cap \PP(1^5,2^4)\) is empty and  also the base locus of the linear system \(|\Oh(2)|\) of quadrics is empty; \(X\) is a smooth Fano 4-fold.
\section{Geography of smooth Fano 4-folds}
\noindent The deformation type of a smooth Fano variety \(X\) depends on  the plurigenera \(h^{0}(-nK_X)\)  of the Hilbert series \(\sum_{n\ge 0}h^0(-nK_X)t^n \) of \(X\), as they are invariant under smooth projective deformations \cite{Inv:pluri}. Therefore if two varieties have different plurigenera then they represent two distinct different deformation families of Fano 4-folds.    Now by using various vanishing theorems one can get the following form of  the Riemann--Roch formula for smooth Fano 4-folds of index 1 \cite[p.48]{Kuchle-97}.
\begin{equation}\chi(-nK_X)=h^0(-nK_X)=1+\dfrac{n(n+1)}{24}(-K_X)^2c_2(X)+\dfrac{n^2(n+1)^2}{24}(-K_X)^4
\end{equation}   
In particular, \begin{equation} h^0(-K_X)=1+\dfrac{(-K_X)^2 c_2(X)}{24}+\dfrac{(-K_X)^4}{6}.\end{equation}
Thus the intersection number \((-K_X)^2c_2(X)\) is  determined easily from the first term if we can compute    \(h^{0}(-K_X)\) and \((-K_X)^4\). In our case, these two invariants can be readily computed from the Hilbert series of \(X\). Therefore we have  the following result. 

\begin{prop}\label{distinguish} Let \(X\) and \(Y\) be two smooth Fano 4-folds of index 1 such that \(h^0(-K_X)\ne h^0(-K_{Y})\) and $(-K_X)^4\ne (-K_{Y})^4$ then \(X\) and \(Y\) belong to two distinct deformation families of smooth Fano 4-folds. 
\end{prop} 
\subsection{Geography with respect to \(h^{0}(-K_X)\) and \((-K_X)^4\)} 
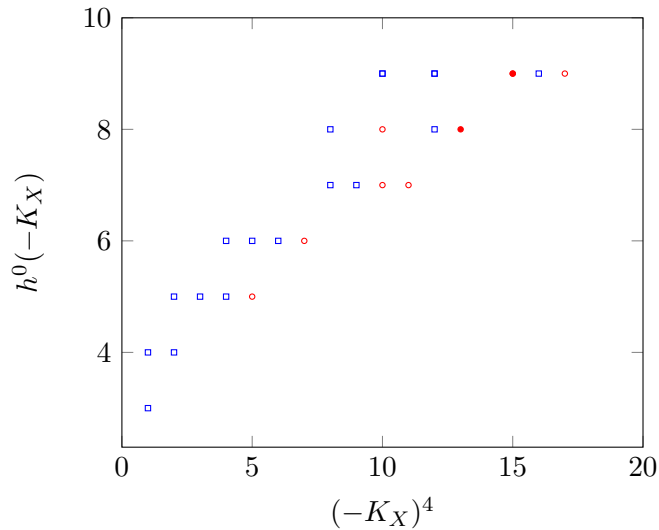
\begin{figure}[htbp]  
  
 \caption{Fano 4-folds with small invariants satisfying \eqref{invariants}. The red circles represent  new examples constructed  in this paper, blue squares represent  the  examples already appeared in \cite{Kuchle-95,Kuchle-97,HLM}, and the   filled red dots represents those  examples which appeared    in this paper and in \cite{Kuchle-97}.   }       \label{diag}
\begin{tikzpicture}
\begin{axis}[xlabel=$(-K_X)^4$,
xmin= 0, xmax= 20, ylabel={$h^0(-K_X)$},
ymax=10]
\addplot[
color = red,
fill = red,
mark size = 1pt,
mark=o,
 % A filled circle
only marks] 
coordinates {

( 10,7 )
( 7,6 )
( 5,5 )
( 15,9 )
( 10,8 )
( 17,9 )
( 11,7 )
};
\addplot[
color = blue,
fill = blue,
mark size = 1pt,
mark=square,
 % A filled circle
only marks] 
coordinates {
( 1, 3 )
( 1, 4 )
( 2, 4 )( 2, 5 )( 3, 5 )( 4, 6 )( 4, 5 )( 5, 6 )( 6, 6 )( 8, 7 )(8,8)( 9, 7 )(10,9)( 12, 8 )(12,9)( 16, 9 )
%Hausen
(12,9)(10,9)
};
\addplot[
color = red,
fill = red,
mark size = 1pt,
% A filled circle
only marks] 
coordinates {
(15,9) (13,8)};

\end{axis}

\end{tikzpicture}
\end{figure}

In total there are at least $987$ known  deformation families of smooth Fano 4-folds with distinct anti-canonical degree  $(-K_X)^4$ and $h^0(-K_X)$.  They  can be found in   \cite{Kuchle-95,Kuchle-97,Q-periods,CKP,CKP-LI,Klash,HLM}. Among these, there are $13$ examples of Fano 4-folds of index 1 which are hypersurfaces or complete intersections in weighted projective spaces \cite{Kuchle-97}. For all these examples, the first plurigenus \(h^0(-K_X)\) and the anti-canonical degree \((-K_X)^4\) satisfy,
\begin{equation}\label{invariants}1\le -K^4\le 17 \textrm{ and } 3\le h^0(-K_X)\le 9.\end{equation} 
  In what follows, we call a smooth  Fano 4-fold to have small invariants if  its invariants satisfy \eqref{invariants}. There are very few families of smooth Fano 4-folds  with such small invariants, other than  hypersurfaces and complete intersections in weighted projective spaces.    In total, excluding the hypersurfaces or complete intersections, only   $7 $ out of the rest of  about $970$ known examples have small invariants. Among these, 3 of them appeared in \cite{Kuchle-95, Kuchle-97}, $3$ are listed in \cite{Q-periods} and $3$ in \cite{HLM} but 2 of them have same invariants as those  in \cite{Kuchle-97}.   All our new families of Fano 4-folds have small invariants and thus they lie    in the left-down corner of the graph of geography of smooth Fano 4-folds if we draw a graph in positive quadrant with \((-K_X)^4\) on \(x\)-axis and \(h^0(-K_X)\) on \(y\)-axis. In figure \ref{diag}, we list those smooth families of Fano 4-folds with known small invariants. 
%%%%%%%%%%%%%%%%%%%%%%%%%%%%%%%%%%%%%%%%%%%%%%%%%%%%%%%%%%
    
  \section*{Acknowledgements} I am thankful to Gavin Brown and Alexander Kasprzyk for helpful discussions.  This research is supported by the  grant No. SR191006 of the   Deanship of Scientific Research (DSR) at King Fahd University of Petroleum and Minerals.  
\bibliographystyle{amsplain}
\bibliography{References}
%--------------------------------------------------------------------------------------------

\end{document}